\documentclass[11pt]{article}
\usepackage{latexsym}
        \usepackage{amscd}
        \usepackage{amssymb}
        \usepackage{amsmath}
        \usepackage{amsthm}
        \usepackage{enumerate}
\def\co{\colon\thinspace}

\newtheorem{thm}{Theorem}[section]

\newtheorem{cor}[thm]{Corollary}

\newtheorem{lem}[thm]{Lemma}

\newtheorem{Example}[thm]{Example}
\newenvironment{ex}{\begin{Example}\rm}{\end{Example}}
\newtheorem{Counterexample}[thm]{Counterexample}

\newtheorem{remark}[thm]{Remark}
\newenvironment{rmk}{\begin{remark}\rm}{\end{remark}}
\newtheorem{Fact}[thm]{Fact}

\newtheorem{Nothing}[thm]{$\!\!\!$}

\newcommand{\be}{\begin{equation}}
\newcommand{\ee}{ \end{equation}} 
\newcommand{\ba}{\begin{eqnarray}}
\newcommand{\ea}{\end{eqnarray}}
\newcommand{\ban}{\begin{eqnarray*}}
\newcommand{\ean}{\end{eqnarray*}}
\newcommand{\pt}{\partial}
\newcommand{\lp}{\langle}
\newcommand{\rp}{\rangle}

\newcommand{\Ric}{\mbox{Ric}}
\newcommand{\Iso}{\mbox{Iso}}

\begin{document}
\abovedisplayskip=6pt plus3pt minus3pt
\belowdisplayskip=6pt plus3pt minus3pt
\title{\bf Metrics of positive Ricci curvature on bundles\rm}
\author{Igor Belegradek \and 
Guofang Wei\thanks {Partially supported by NSF Grant \# DMS-9971833.}}
\date{}
\maketitle
\begin{abstract}
We construct new examples of manifolds of positive Ricci 
curvature which, topologically, are vector bundles over compact
manifolds of almost nonnegative Ricci curvature. 
In particular, we prove that if $E$ is the total space of a vector bundle 
over a compact manifold of nonnegative Ricci curvature,
then $E\times\mathbb R^p$ admits a complete metric of positive Ricci 
curvature for all large $p$. 
\end{abstract}

\section{Introduction}

According to the soul theorem of J.~Cheeger and D.~Gromoll, 
a complete open manifold of nonnegative sectional curvature,
denoted $K\ge 0$, is the total space of a vector bundle over 
a compact manifold with $K\ge 0$. 
Manifolds of nonnegative Ricci curvature are 
much more flexible, and nowadays 
there are many examples of complete manifolds 
of $\Ric\ge 0$  which are not even homotopy equivalent
to complete manifolds of $K\ge 0$.
These include manifolds not homotopy equivalent to closed 
manifolds~\cite{GM}, 
manifolds not satisfying Gromov's Betti numbers estimate~\cite{SY},
manifolds of infinite topological type 
(see~\cite{Men} and references therein),
manifolds with not virtually-abelian fundamental group~\cite{wei},
compact spin Ricci-flat {$4$-manifolds} 
with nonzero $\hat A$-genus~\cite[6.27]{Bes},~\cite{lot},
and complements to certain smooth divisors in compact 
K\"ahler manifolds~\cite{TY}.

In~\cite{BW} the authors constructed 
the first (to our knowledge) examples of 
complete manifolds with $\Ric>0$ which
are homotopy equivalent but not homeomorphic to complete 
manifolds with $K\ge 0$. Topologically, the manifolds are
vector bundles over tori (in~\cite{BW}
we also constructed vector bundles
over nilmanifolds carrying $\Ric>0$). 
According to~\cite{OW, BK2, BK3}
in each rank only finitely many vector bundles over tori admit metrics
with $K\ge 0$, and more generally, a majority of vector bundles
over a fixed manifold $B$ 
with $K\ge 0$ does not admit a metric with $K\ge 0$,
provided $B$ has a sufficiently large first Betti number.

In this paper we greatly extend the results of~\cite{BW}. 
In particular, in~\cite{BW} we asked whether most vector 
bundles over compact manifolds with $\Ric\geq 0$ admit 
metrics with $\Ric >0$, and in Corollary~\ref{cor: vect over ric nneg} 
we show that this is true stably, i.e. after multiplying by some
high-dimensional Euclidean space.

Let $\mathcal B$ be the smallest class of manifolds 
containing all compact manifolds with $\Ric\ge 0$, and 
any manifold $E$ which is the total space of a smooth
fiber bundle $F\to E\to B$ where $B\in\mathcal B$,
$F$ is a compact manifold of nonnegative Ricci curvature, 
and the structure group of the bundle
lies in the isometry group of $F$. 
In other words, a manifold in $\mathcal B$ is an 
iterated fiber bundle
such that all the fibers, and the base at the very first step 
are compact manifolds of $\Ric\ge 0$, and the structure groups
lie in the isometry groups of the fibers. Note that $\mathcal B$
is closed under products.
Here is our main result. 

\begin{thm}
\label{thm: vb over b} Let $B\in\mathcal B$, and
let $E(\xi)$ be the total space of a vector bundle $\xi$ over $B$.
Then $E(\xi)\times \mathbb R^p$ admits a complete Riemannian metric of
positive Ricci curvature for all sufficiently large $p$.
\end{thm}

\begin{cor}\label{cor: vect over ric nneg}
Let $B$ be a compact manifold with $\Ric\ge 0$. 
If $E(\xi)$ is the total space of a vector bundle $\xi$
over $B$, then $E(\xi)\times \mathbb R^p$ 
admits a complete Riemannian metric of
positive Ricci curvature for all large $p$.
\end{cor}

The case when $\xi$ is a trivial vector bundle in Theorem~\ref{thm: vb over b}   
generalizes the main result of~\cite{wei}. 
Actually, in Theorem~\ref{thm: vb} and Section~\ref{sec: ibii}
we prove a version of Theorem~\ref{thm: vb over b} 
for a larger class of base manifolds 
including some iterated fiber bundles 
with almost nonnegatively curved fibers.
(A special case of this result with nilmanifolds as fibers
was proved in~\cite[Theorem 4]{W2}).

Theorem~\ref{thm: vb over b} should be compared with the results
of J.~Nash and L.~Berard-Bergery~\cite{Nash, BBer} who proved that
the class of compact manifolds with $\Ric>0$ is closed under taking fiber
bundles with structure groups lying in the isometry groups of the fibers,
i.e. if $B,F$ admit $\Ric>0$, then so does $E$.
Furthermore, 
any vector bundle of rank $\ge 2$ over a compact manifold with $\Ric>0$ 
carries a complete metric with $\Ric>0$~\cite{Nash, BBer}. Note that
rank one vector bundles cannot carry $\Ric>0$
by the Cheeger-Gromoll splitting theorem.

By considering iterated bundles such that the fibers
have large isometry groups
(e.g. if the fibers are spheres, or compact Lie groups), 
one sees that $\mathcal B$ 
contains many different topological types.
For example, $\mathcal B$ contains all nilmanifolds,
or more generally, all iterated linear sphere bundles,
as well as all iterated principal bundles (with compact fibers)
over compact manifolds with $\Ric\ge 0$.

All manifolds in $\mathcal B$ admit metrics of almost
nonnegative Ricci curvature (which was certainly known to Nash
and Berard-Bergery),
however, many of the manifolds in $\mathcal B$ 
do not admit metrics with $\Ric\ge 0$ (e.g.~nilmanifolds). 
More examples are given in Section~\ref{sec :fb and ric}
where we show that if $E$ is a compact manifold with $\Ric(E)\ge 0$
which fibers over a torus $T$, then 
the pullback of the bundle $E\to T$ to a finite cover of $T$ 
has a section. 

Theorem~\ref{thm: vb over b} becomes false
if in the definition of $\mathcal B$ we do not assume
that the structure group lies in the isometry group
of the fiber. For example, Theorem~\ref{thm: vb over b}
fails when $B$ is a compact Sol $3$-manifold (which is a
$2$-torus bundle over a circle) because then $\pi_1(B)$
is not virtually nilpotent.

To prove~\ref{thm: vb over b} we show that any manifold in 
$\mathcal B$ admits what we call a metric of 
{\it almost nonnegative Ricci curvature 
with good local basis} (see Section 7 for a precise
definition), and then we prove in Theorem~\ref{thm: positive}
that the product of
any manifold carrying such a metric with a high-dimensional 
Euclidean space admits a complete metric with $\Ric>0$.

As was suggested in~\cite{W2}, it may well be true in general that 
if $E$ is a manifold of almost nonnegative
Ricci curvature, then $E\times\mathbb R^p$ has a complete 
metric with $\Ric>0$ for large $p$. 
Our results is a further step in this direction. 

The minimal value of $p$ coming from our construction generally 
depends on $\xi$, and typically is very large. It would be
interesting to find obstructions for small $p$'s. 
For manifolds with infinite fundamental group there are 
obstructions due to M.~Anderson~\cite{And}, e.g. 
no $\mathbb R^2$-bundle over
a torus admits a complete metric with $\Ric>0$.
By contrast, no obstructions are known in the simply-connected
case, say it is unclear whether the product of $\mathbb R^2$ and
a simply-connected compact Ricci-flat manifold can admit metrics
with $\Ric>0$.

The structure of the paper is as follows.
In sections~\ref{sec: fb and subm}--\ref{sec: scaling}
we review a well-known construction of metrics
of almost nonnegative Ricci curvature on fiber bundles. 
We frequently refer to~\cite[Chapter 9]{Bes} for details.
Sections~\ref{sec: warp1}--\ref{sec: warp2}
contain a curvature computation generalizing
the main computation in~\cite{wei}.
Section~\ref{sec: metrics on fb}
provides an easy route to the proof 
of Corollary~\ref{cor: vect over ric nneg}.
Main technical results are proved in 
sections~\ref{sec: ibi}--\ref{sec: ibii}.
In section~\ref{sec :fb and ric} we give examples
of iterated fiber bundles in the class $\mathcal B$
which do not admit metrics of $\Ric\ge 0$.

It is a pleasure to thank Vitali Kapovitch for help with
Example~\ref{ex: vb error term}, and McKenzie Wang for 
providing a list of misprints in~\cite[Chapter 9]{Bes}.

\section{Fiber bundles and Riemannian submersions}
\label{sec: fb and subm}
In sections~\ref{sec: fb and subm}--\ref{sec: scaling} we
let $\pi\co E\to B$ be a Riemannian submersion of complete 
Riemannian manifolds with totally geodesic fibers.
By~\cite[9.42]{Bes}, $\pi$ is a smooth fiber 
bundle whose structure group 
is a subgroup of the isometry group of the fiber;
all fibers are isometric,
and we denote a typical fiber by $F$~\cite[9.56]{Bes}.

\begin{ex} \label{ex: local metric}
Let $B$, $F$ be complete Riemannian manifolds 
where the metric on $F$ is invariant under a compact Lie group $G$. 
Let $\pi\co E\to B$ be a smooth fiber bundle with fiber $F$ 
and structure group $G$. 
Then there exists a complete Riemannian metric on $E$, 
making $\pi$ a Riemannian submersion with
totally geodesic fibers isometric to $F$~\cite[9.59]{Bes}. 
(To construct the metric, 
think of $E$ as an $F$-bundle associated with a principal
$G$-bundle $P$ over $M$. Choose a connection on $P$ which
defines a horizontal and vertical distributions on $E$.
Introduce the metric on $E$ by making it equal to the metric on 
$B$, $F$ in the horizontal, vertical subspaces, respectively,
and making the subspaces orthogonal).
In particular, setting $F=G$ equipped with a left-invariant metric,
we get a Riemannian submersion $P\to B$ where $P$ is 
any principal $G$-bundle $P$ over $M$. 
\end{ex}

\begin{ex}\label{ex: global metric}
There is another useful metric on $E$. Namely, we think of $E$ as
the quotient $(P\times F)/G$, take the product metric
on $P\times F$ (where the metric on $P$ is defined 
as in Example~\ref{ex: local metric}), and give $E$ the Riemannian submersion 
metric 
(with horizontal spaces orthogonal to the fibers).
The fibers in this metric are also totally geodesic~\cite[3.1]{Nash}
and isometric to the base of the Riemannian submersion 
$(G\times F)\to (G\times F)/G$ which is diffeomorphic to $F$. 
For example, the total space of any rank $n$ vector bundle
gets the Riemannian submersion metric from
$(P\times\mathbb R^n)/O(n)$ where $\mathbb R^n$
is the Euclidean space; the fibers are isometric to 
$(O(n)\times\mathbb R^n)/O(n)$. 
\end{ex}

\begin{ex}\label{ex: flat}
One particularly simple kind of Riemannian submersion is 
a locally isometric flat bundle which can be described
as follows.  
Start with complete Riemannian manifolds $B$, $F$, and 
a homomorphism $\rho\co\pi_1(B)\to\Iso(F)$, and consider
the flat $F$-bundle $E$ over $B$ with holonomy $\rho$.
If $\tilde B$ is the universal cover of $B$, then the 
total space $E$ of the bundle is the quotient of
$\tilde B\times F$ by the $\pi_1(B)$-action 
given by $\gamma(\tilde b,f)=(\gamma(\tilde b), \rho(\gamma)(f))$.
The bundle projection $E\to B$ is induced by the projection
$\tilde B\times F\to B$. The product metric on $\tilde B\times F$
defines a metric on $E$ which makes the projection $E\to B$
into a Riemannian submersion with totally geodesic fibers
isometric to $F$, and zero $A$-tensor.
\end{ex}

\section{Choosing a basis}
\label{sec: basis}
We use specific local trivializations of $\pi$
defined by taking a small strictly convex open ball 
$U$ in $B$ and considering horizontal lifts of the radial
geodesics emanating from the center of $U$. 
Thus $\pi^{-1}(U)$ gets identified with $U\times F$
where the map $\{u\}\times F\to\{u^\prime\}\times F$
given by $(u,f)\to (u^\prime,f)$ is an isometry
with respect to the induced metric on the fibers (see~\cite[9.56]{Bes}).

To simplify the curvature computation we choose
a basis on $E$ as follows. Fix an arbitrary  
point $e$ of $E$. Let $F_e$ be the fiber passing through $e$. 
At the tangent spaces $T_e F_e$, $T_{\pi(e)}B$, 
start with arbitrary orthonormal bases 
$\{\hat W_i\}$, $\{\check H_j\}$, and 
extend $\hat W_i$, $\check H_j$ to orthonormal vector fields 
also denoted $\hat W_i$, $\check H_j$ on neighborhoods of 
$e\in F_e$, $\pi(e)\in B$, respectively,
such that $\lp [\check H_i,\check H_j], \check H_i\rp$ and
$\lp [\hat W_i,\hat W_j], \hat W_i\rp$ vanish at $e$.
(This can be achieved, for example, by choosing
the extension so that $[\hat W_i,\hat W_j]$,
$[\check H_i,\check H_j]$ vanish at $e$ or geodesic frame. 
The reason we care about the property
can be seen in Lemma~\ref{lem: dr terms}).

Use the above local trivializations to extend these
vector fields to vector fields $W_i$, $H_j$ defined 
on a neighborhood of $e$ in $E$. 
Thus, at any point of the neighborhood, $W_i$ is vertical, and 
$H_j$ is horizontal with $\pi_*(H_j)=\check{H_j}$.
We conclude that $\{W_i,H_j\}$ is an orthonormal basis 
on a neighborhood of $e\in E$. 

This basis has the property that 
$\lp [A,B],A\rp\vert_e=0$ for any $A,B\in\{W_i,H_j\}$.
Indeed, at $e$ the following is true. 
First, since $\pi_*[H_i,H_j]=[\check H_i,\check H_j]$,
we get $\lp [H_i,H_j], H_i\rp =
\lp [\check H_i,\check H_j], \check H_i\rp=0$.
Similarly, $\pi_*W_j=0$ implies  $\lp [H_i,W_j], H_i\rp =0$.
By construction $\lp [W_i,W_j], W_i\rp= 
\lp [\hat W_i,\hat W_j], \hat W_i\rp=0$.
Finally, by the Koszul's formula
\[\lp [W_i,H_j], W_i\rp=-\lp\nabla^E_{W_i}W_i,H_j\rp=0\] because 
the fibers are totally geodesic hence $\nabla^E_{W_i}W_i$ is vertical.

\section{Submersions with bounded error term}

We say that a Riemannian submersion with totally geodesic fibers
has a {\it bounded error term} if there is an (independent of a point)
constant $C>0$ such that the tensors
$\lp AU,AV\rp$, $\lp A_{X},A_{Y}\rp$, 
$\lp\check\delta A(X),U\rp$, defined in~\cite[9.33]{Bes},
are bounded by $C$ in absolute value, where 
$U, V, X, Y$ are unit vector fields, $U, V$ are vertical, and
$X, Y$ are horizontal. 
If $E$ is compact, then by continuity $\pi$ has a bounded
error term, but there are some other examples.
In fact, by~\cite[9.36]{Bes}, we have
\begin{eqnarray*}
& \lp AU,AV\rp=
\Ric_E(U,V)-\Ric_F(\hat U,\hat V), \\
& 2\lp A_{X},A_{Y}\rp=
\Ric_B(\check X,\check Y)-\Ric_E(X,Y), \\ 
& \lp\check\delta A(X),U\rp=-\Ric_E(X, U),
\end{eqnarray*} 
so if $\Ric_B$, $\Ric_F$, $\Ric_E$ are bounded
in absolute value, then $\pi$ has a bounded error term
(Here $\check X$ denotes $\pi_* X$, and $\hat U$
denotes the restriction of $U$ to the fiber, etc).

\begin{ex}\label{ex: vb error term}
The projection of any rank $n$ vector bundle $E\to B$ over a 
compact manifold $B$ can be made into a Riemannian submersion 
as in Example~\ref{ex: global metric}. 
Here $E$ gets the metric as the base of the Riemannian submersion
$p\co (P\times\mathbb R^n)\to (P\times\mathbb R^n)/O(n)$. 
By a simple computation (done e.g. in \cite[Page 361]{A2}), the $A$-tensor of 
$p$ is such that
$|A_X Y|$ is uniformly bounded above
for any $p$-horizontal orthogonal unit vector fields $X,Y$.
Then since the sectional curvature of $P\times\mathbb R^n$ is bounded
above and below, so is the sectional curvature of $E$  
by the O'Neill's formula. By the same argument the sectional curvature
of $(O(n)\times\mathbb R^n)/O(n)$, which is the fiber of $E\to B$, 
is bounded above and below, hence the submersion
$E\to B$ has a bounded error term.
\end{ex}

\section{Scaling Riemannian submersions}
\label{sec: scaling}
Given $t\in (0,1]$, let $g_{t,E}$ be the family of metrics on $E$ 
defined by
\begin{equation}
g_{t,E}(U+X,V+Y)=t^2g_E(U,V)+g_E(X,Y),
\end{equation}
where $U,V$ are vertical, and $X,Y$ are horizontal.
We refer to~\cite[9G]{Bes} for more information on such metrics
(note that Besse scales the metric on fibers by $t$ while
we use $t^2$).
By~\cite[9.70]{Bes} (note the wrong sign in 9.70c), 
applied to the submersion 
$\pi\co E\to B$, we get 
\begin{eqnarray*}
& \Ric_{g_{t,E}}(U,V)=\Ric_{g_{t,F}}(\hat U,\hat V)+t^4\lp AU,AV\rp, \\
& \Ric_{g_{t,E}}(X,Y)=\Ric_{g_{B}}(\check X,\check Y)-
2t^2\lp A_{X},A_{Y}\rp ,\\
& \Ric_{g_{t,E}}(X,U)=-t^2\lp \check\delta A(X),U\rp .
\end{eqnarray*}
Here $g_{t,F}$ is the metric on $F$ induced by $g_{t,E}$,
in particular, $\Ric_{g_{t,F}}=t^{-2}\Ric_F$.

For example, if in the construction of $\{H_i, W_j\}$, 
we assume that $\{\hat W_i\}$, $\{\check H_j\}$ are 
orthonormal bases of eigenvectors diagonalizing 
the Ricci tensors $Ric_{F_e}$, $Ric_B$
(so that for $i\neq j$, 
$Ric_F(\hat W_i,\hat W_j)=0=Ric_B(\check H_i,\check H_j)$),
and if $\pi$ has a bounded error term, then
\begin{eqnarray}
&|\Ric_{g_{t,E}}(A,B)|\le Ct \ \mbox{for any distinct}\ A,B\in \{H_i, 
\frac{W_j}{t}\}, \label{eqn: 5.2}\\
& \Ric_{g_{t,E}}(\frac{W_i}{t},\frac{W_i}{t})\ge \Ric_F(\hat W_i,\hat W_i)t^{-2}, 
\label{eqn: 5.3}\\
& \Ric_{g_{t,E}}(H_i,H_i)\ge \Ric_B(\check H_i,\check H_i)-Ct^2 \label{eqn: 5.4}.
\end{eqnarray}
So if $\Ric(F)\ge 0$ and $\Ric(B)\ge 0$,
then $\Ric_{g_{t,E}}(W_i,W_i)\ge 0$ and $\Ric_{g_{t,E}}(H_i,H_i)\ge -Ct^2$.

\section{Computing the Ricci curvature of a warped product}
\label{sec: warp1}
In this section we present a generalization of
the main computation of~\cite{wei}.
Let $E$ be an $n$-manifold with a smooth family of complete Riemannian
metrics $g_r$, $r\ge 0$. 
Let $f$ be a function of $r$ to be specified later.
Consider a metric on $N=E\times S^{p-1}\times (0,\infty)$
given by $g=s_r+dr^2$ where $s_r=g_r+f^2ds^2$ and 
$ds^2$, $dr^2$ are the canonical
metrics on $S^{p-1}$, and $(0,\infty)$. 

Note that $N$ is an open subset in a manifold diffeomorphic to 
$E\times\mathbb R^p$.
Our goal is to find the conditions ensuring 
that $g$ extends to a complete metric 
on $E\times\mathbb R^p$ with $Ric(g)>0$ for all 
sufficiently large $p$.

For each point $e\in E$ assume that
there is a basis of vector fields $\{X_i\}$
on a neighborhood $U_e\subset E$ of $e$ such that 
for any $r\ge 0$, $g_r(X_i,X_j)=0$ on $U_e$ if $i\neq j$.
Fix one such a basis for each $e$.
Let $h_i(r)=\sqrt{g_r(X_i,X_i)}$ so that $Y_i=X_i/h_i$
form a $g_r$-orthonormal basis on $U_e$ for any $r\ge 0$.
Since $X_i\neq 0$ and $g_r$ is nondegenerate, $h_i$ 
is a positive function on $[0,\infty)$.
Assume furthermore that $h_i$ is smooth and 
$h_i^{(odd)}(0)=0$ (i.e. all odd derivatives of $h_i$
at zero vanish).
An argument similar to~\cite[Page 13]{Pet} 
shows that $g$ is a smooth complete Riemannian 
metric on $E\times\mathbb R^p$ if $f(0)=0$, $f(r)>0$ for $r>0$, 
$f^\prime(0)=1$, $f^{(even)}(0) = 0$.

For the rest of the section we 
fix an arbitrary point $(e, s, r)\in N$. 
Our goal is to compute the Ricci tensor of $g$ 
at the point in terms of $f$, $h_i$'s, and the Ricci 
tensor of $g_r$.

At a neighborhood of $s\in (S^{p-1},ds^2)$ choose 
an orthonormal frame $\{V_j\}$ 
with $\nabla^{ds^2}_{V_j}V_k\vert_{s}=0$ for any $j,k$
where $\nabla^{ds^2}$ is the Riemannian connection for $ds^2$.
In particular, $[V_i,V_j]\vert_s=0$.
Denote $U_j=V_j/f$ and $\pt_r=\frac{\pt}{\pt r}$ and so that 
\[\beta_r=\{\pt_r, U_1,\dots U_{p-1}, Y_1, \dots, Y_n\},\qquad
\beta=\{U_1,\dots U_{p-1}, Y_1, \dots, Y_n\}\]
are orthonormal frames at $(e,s,r)$ and $(e,s)$, respectively.

We start with a few elementary observations which we use throughout this
section. Any function of $r$ has zero derivative in the direction of
$X_i$ (because $X_i$ is independent of $r$), or $V_i$ (since $V_i$
is tangent to the level surface of the function).
Hence $[Y_i,Y_j]=\frac{1}{h_ih_j}[X_i,X_j]$ and 
$[U_i,U_j]\vert_s=\frac{1}{f^2}[V_i,V_j]\vert_s=0$.
Also the flow of $\pt_r$ preserves $X_i$, $V_j$
so $[\pt_r,X_i]=0=[\pt_r, V_j]$ which implies 
$[Y_i,\pt_r]=\frac{h_i^\prime}{h_i}Y_i$ and
$[U_j,\pt_r]=\frac{f^\prime}{f}U_j$. 

We frequently use the Koszul's formula for the Riemannian connection of $g$
computed in the orthonormal basis, namely for any $A,B,C\in\beta$
\[2\lp\nabla_A B,C\rp_g=\lp[A,B],C\rp_g+\lp[C,A],B\rp_g+\lp[C,B],A\rp_g.\]
In particular, the Koszul's formula implies
$\nabla_{Y_i}\pt_r=\frac{h_i^\prime}{h_i}Y_i$
and $\nabla_{U_j}\pt_r=\frac{f^\prime}{f}U_j$
so that $\nabla_{\pt_r} {Y_i}$, $\nabla_{\pt_r} U_j$ vanish.
Also $\nabla_{\pt_r}\pt_r=0$. 
For any $i,j$, $[U_i,U_j], [Y_i, U_j]$,
$[Y_i,Y_j]$ are all tangent to $E\times S^{p-1}$
so they all have zero $\pt_r$-component.
Therefore, by the Koszul's formula 
the terms $\lp\nabla_{Y_i}Y_j, \pt_r\rp$, $\lp\nabla_{Y_i}U_j,\pt_r\rp$, 
$\lp\nabla_{U_i}U_j, \pt_r\rp$ all vanish for $i\neq j$.
The Koszul's formula also implies
$\lp\nabla_{Y_i}Y_i,\pt_r\rp=\lp[\pt_r,Y_i],Y_i\rp=-\frac{h_i^\prime}{h_i}$,
and similarly $\lp\nabla_{U_i}U_i,\pt_r\rp=-\frac{f^\prime}{f}$.

The second fundamental form of the submanifold $E\times S^{p-1}\subset N$
is given by $II(A,B)=\lp\nabla_A B,\pt_r\rp \pt_r$ for $A,B\in\beta$. 
Thus, $II(Y_i,Y_i)=-\frac{h_i^\prime}{h_i}\pt_r$,
$II(U_i,U_i)=-\frac{f^\prime}{f}$ while all the mixed terms 
$II(Y_i,U_j)$ as well as the terms
$II(U_i,U_j)$, $II(Y_i,Y_j)$ for $i\neq j$ vanish.
Now the Gauss equation~\cite[VII.4.1]{KN} gives the sectional
curvature of $N$ in terms of the sectional curvature of 
$(E\times S^{p-1}, s_r)$ as follows:
\begin{eqnarray*}
 & K_g(Y_i,U_j)=K_{s_r}(Y_i,U_j)-\frac{h_i^\prime f^\prime}{h_i f}\\
 & K_g(Y_i,Y_j)=K_{s_r}(Y_i,Y_j)-\frac{h_i^\prime h_j^\prime}{h_ih_j}\\
 & K_g(U_i,U_j)=K_{s_r}(U_i,U_j)-\left(\frac{f^\prime}{f}\right)^2. 
\end{eqnarray*}
Since $s_r$ is the product metric on $E\times S^{p-1}$, we
get $\nabla^{s_r}_{U_i}Y_j=0=\nabla^{s_r}_{Y_j}U_i$ so that
the curvature $(2,2)$-tensor $R_{s_r}(U_i,Y_j)$ vanishes.
Hence $K_{s_r}(U_i,Y_j)=0$ and $\lp R_{s_r}(A,B)C,D\rp=0$
for any $A,B,C,D\in\beta$ unless $A,B,C,D$ are all tangent
to $E$ or $S^{p-1}$. Also 
\begin{eqnarray*}
 \lp R_{s_r}(Y_i,Y_k)Y_l,Y_j\rp=\lp R_{g_r}(Y_i,Y_k)Y_l,Y_j\rp,
& \lp R_{s_r}(U_i,U_k)U_l,U_j\rp=\lp R_{f^2ds^2}(U_i,U_k)U_l,U_j\rp,\\
 K_{s_r}(Y_i,Y_j)=K_{g_r}(Y_i,Y_j),
& K_{s_r}(U_i,U_j)=K_{f^2ds^2}(U_i,U_j)=1/f^2. 
\end{eqnarray*}

Now we turn to the terms of the curvature tensor involving 
$\pt_r$ and the mixed terms.
\begin{lem}
(1) $R_g(\pt_r, U_i)\pt_r=\left(\frac{f^{\prime\prime}}{f}\right) U_i$ and 
$R_g(\pt_r, Y_i)\pt_r=\left(\frac{h_i^{\prime\prime}}{h_i}\right) Y_i$;\\
(2) $K_g(U_i,\pt_r)=-\left(\frac{f^{\prime\prime}}{f}\right)$ and
$K_g(Y_i,\pt_r)=-\left(\frac{h_i^{\prime\prime}}{h_i}\right)$;\\
(3) $\Ric_g(Y_i,Y_j)=\Ric_{g_r}(Y_i,Y_j)$ 
and $\Ric_g(U_i,U_j)=0$ for $i\neq j$;\\
(4) $\Ric_g(U_i,Y_j)=0$.
\end{lem}
\begin{proof}
(1) follows directly from $\nabla_{\pt_r}\pt_r=0$, 
$\nabla_{Y_i}\pt_r=\frac{h_i^\prime}{h_i}Y_i$, 
$\nabla_{U_j}\pt_r=\frac{f^\prime}{f}U_j$,
$[Y_i,\pt_r]=\frac{h_i^\prime}{h_i}Y_i$, and
$[U_j,\pt_r]=\frac{f^\prime}{f}U_j$.

(2) follows from (1). 

(3)-(4) By the Gauss equation~\cite[VII.4.1]{KN}, we have
$R_g(A,C,C,B)=R_{s_r}(A,C,C,B)$ for any $A,B,C\in\beta$, $A\neq B$.
By part (1), $R_g(A,\pt_r)\pt_r$ is proportional to $A$, so
$R_g(A,\pt_r,\pt_r, B)=0$ if $A\neq B$. Now $\Ric_g(A,B)$ is the sum
of $R_g(A,\pt_r,\pt_r, B)$ and all the terms $R_g(A,C,C,B)$ so we get 
$\Ric_g(A,B)=\Ric_{s_r}(A,B)$. 
Using the above formulas for $R_{s_r}$, with $i\neq j$, 
we get $\Ric_{s_r}(Y_i,Y_j)=\Ric_{g_r}(Y_i,Y_j)$, and
$\Ric_{s_r}(U_i,U_j)=\Ric_{f^2ds^2}(U_i,U_j)=0$
where the last equality is true since $f^2ds^2$ has constant sectional
curvature. Finally, $\Ric_{s_r}(Y_i,U_j)=0$ since
all the involved curvature tensor terms vanish. 
\end{proof}

{\bf Correction} (added on August 28, 2010): Parts (1) and (4) of 
Lemma~\ref{lem: dr terms}, and formula (\ref{form: dr ric})
below are incorrect. A correction can be found in
Appendix C of~\cite{Bel-rh-warp} where it is explained
why the mistake does not affect other results
of the present paper.

\begin{lem} \label{lem: dr terms}
(1) $\lp R_g(\pt_r,Y_i)Y_i,Y_j\rp=\lp [Y_i,Y_j],Y_i\rp
\left(\frac{h_i^\prime}{h_i}+\frac{h_j^\prime}{h_j}\right)$;\\
(2) $\lp R_g(\pt_r,A)B,C)=0$ for any $A, B, C\in\beta$
unless $A, B, C$ are all tangent to $E$;\\
(3) $\Ric_g(\pt_r,U_i)=0$;\\
(4) $\Ric_g(\pt_r,Y_j)=\sum_{i\neq j}\lp [Y_i,Y_j],Y_i\rp
\left(\frac{h_i^\prime}{h_i}+\frac{h_j^\prime}{h_j}\right)$.
\end{lem}
\begin{proof} Using the Koszul's formula, we compute 
the components of $\nabla_{Y_i}Y_i$ in the basis $\beta_r$
to deduce that 
$\nabla_{Y_i}Y_i=\sum_{j}\lp [Y_j,Y_i],Y_i\rp Y_j-\frac{h_i^\prime}{h_i}\pt_r$.
Then using that $\nabla_{\pt_r}Y_i=0$, 
$[Y_i,\pt_r]=\frac{h_i^\prime}{h_i}Y_i$, we can compute the
curvature tensor $R_g(\pt_r,Y_i)Y_i$, and then its components in $\beta$
so that (1)-(2) follow by a straightforward computation. 
A similar argument gives the result for $U_i$'s the only difference
being $[U_i,U_j]=0$ at $(e,s,r)$. Finally, (3)-(4) follow from
(1)-(2) by substitution.
\end{proof}

In summary, we have
\begin{equation}\label{eq: ric-first}
\Ric_g(U_i,U_i)=(p-2)\frac{1-(f^\prime)^2}{f^2}-\frac{f^\prime}{f}
\sum_{i=1}^{n}\frac{h_i^\prime}{h_i}-\frac{f^{\prime\prime}}{f}
\end{equation}

\begin{equation}
\Ric_g(Y_i,Y_i)=\Ric_{g_r}(Y_i,Y_i)-
(p-1)\frac{f^\prime h_i^\prime}{fh_i}-\frac{h_i^\prime}{h_i}
\sum_{k\neq i}\frac{h_k^\prime}{h_k}-\frac{h_i^{\prime\prime}}{h_i}
\end{equation}
\begin{equation}
\Ric_g(\pt_r,\pt_r)=
-(p-1)\frac{f^{\prime\prime}}{f}-
\sum_{i=1}^{n}\frac{h_i^{\prime\prime}}{h_i}
\end{equation}
\begin{equation}
\Ric_g(U_i,U_j)=0 \ \mbox{for} \  i\neq j, \ \ 
\Ric_g(\pt_r,U_i)=\Ric_g(U_i,Y_j)=0 \  \mbox{for all} \ i,j,
\end{equation}

\begin{equation}\label{form: dr ric}
\Ric_g(\pt_r,Y_j)=\sum_{i\neq j}\lp [Y_i,Y_j],Y_i\rp
\left(\frac{h_i^\prime}{h_i}+\frac{h_j^\prime}{h_j}\right),
\end{equation}
\begin{equation}\label{eq: ric-last}
\Ric_g(Y_i,Y_j)=\Ric_{g_r}(Y_i,Y_j) \ \ \mbox{for} \  i\neq j.
\end{equation}

\section{Turning almost nonnegative into positive}
\label{sec: warp2}

For $q>0$ and $c,m\ge 0$, define the class 
$\mathcal M_{q}(c,m)$ of smooth manifolds by requiring that 
any $E\in\mathcal M_q(c,m)$ has a smooth family
of complete Riemannian metric $g_t$, $t\in (0,1]$ such that any $e\in E$
has a neighborhood $U$ with a basis of vector fields $X_i$
such that \\
(1) $Y_i=X_i/t^{m_i}$ form a $g_t$-orthonormal basis
on $U$ for some $m_i\in [0,m]$, depending on $e$,\\
(2) at the point $e$ the following holds: 
$\lp[Y_i,Y_j],Y_i\rp_{g_t}=0$,
$\Ric_{g_t}(Y_i,Y_i)\ge -ct^q$, for all $i,j$, and  
$|\Ric_{g_t}(Y_i,Y_j)|\le ct^q$ for $i\neq j$.

Let $\mathcal M_{q} = \cup_{c,m}\mathcal M_{q}(c,m)$ 
and $\mathcal M = \cup_q \mathcal M_{q}$. 
We refer to $\mathcal M$ as the class of 
{\it almost nonnegatively Ricci curved manifolds with good 
local basis}. We show in Theorem~\ref{thm: positive} that if 
$E\in\mathcal M$, then $E\times\mathbb R^p$ admits 
a complete metric with $\Ric>0$ for all large $p$. 

Reparametrizing $t$ by $t^r$ for $r>0$, we get 
$\mathcal M_{rq}(c,m) = \mathcal M_q(c,rm)$. 
Also rescaling $g_t$ by $t^r$, we have 
$\mathcal M_{q-2r}(c,m) = \mathcal M_q(c,r+m)$ for $0<r< q/2$. 
Note that $\mathcal M_q(c,m)\subset\mathcal M_s(c,m)$ if $q\ge s$. 

\begin{ex} \label{ex: ric nneg}\ \\
(1) If $E$ has a complete metric $g$ with $\Ric\ge 0$, then 
$E\in\mathcal M_q(0,0)$ for every $q$.
(Take $g_t=g$ for all $t$, 
choose a $g$-orthonormal basis $\{X_i\}$
of $T_eE$ diagonalizing
the Ricci tensor at $e$, and extend it to an orthonormal 
basis of vector fields on $U$ satisfying $[X_i,X_j]=0$).\\
(2) According to~\cite{W2}, any $n$-dimensional nilmanifold,
(or even any $n$-dimensional infranilmanifold~\cite[Theorem 3]{W2}) 
lies in $\mathcal M_q(c,2^{n-2}(q-1) +1)$ 
for every $q$, where $c$ depends on the 
structure constants and $n$.
\end{ex}

Using the formulas of Section~\ref{sec: warp1}, we now 
prove the following.
\begin{thm}  \label{thm: positive}
For each $E\in\mathcal M$ there is a complete metric on
 $E\times\mathbb R^p$ 
 with $\Ric>0$ for all $p$ large. 
 \end{thm}

If $E\in\mathcal M$, then $E\in\mathcal M_q$ for some $q>0$. 
By reparametrizing 
and rescaling we can assume that $E\in\mathcal M_2$ with all $m_i$ 
positive. Now Theorem~\ref{thm: positive}
follows from the theorem below.  

\begin{thm} \label{thm: specific warp} 
If $E \in \mathcal M_2(c,m)$ is an $n$-manifold
with all $m_i$ positive (equivalently, $E \in \mathcal M_q$, where $q >2$), 
then there is an explicit function 
$k(n,c,m)$ such that for all $p>k(n,c,m)$,
there is a complete Riemannian metric $g$ on $E\times\mathbb R^p$
of positive Ricci curvature.
\end{thm}

\begin{proof}
Since $E \in \mathcal M_2(c,m)$ with all $m_i$ positive, by parametrizing 
$t=h(r)$ where $h(r)=(1+r^2)^{-1}$, $r>0$, 
we get a smooth family of complete Riemannian 
metrics $g_r$, $r>0$, such that
at a neighborhood $U$ of any 
point of $e\in E$ there is a basis of vector fields $\{X_i\}$
such that $Y_i=X_i/h^{m_i}$ form a $g_r$-orthonormal
basis on $U$ where $m_i\in (0,m]$,
and at $e$ we have $\lp [Y_i,Y_j],Y_i\rp_{g_r}=0$ for all $i,j$, and  
$|\Ric_{g_r}(Y_i,Y_j)|\le c h^2$ if $i\neq j$, and 
$\Ric_{g_r}(Y_i,Y_i)\ge -c h^2$ for all $i$.

Let $g=g_r+ dr^2+ f^2ds^2$ be the metric on $E\times\mathbb R^p$, where 
$ds^2$ are the canonical
metrics on $S^{p-1}$ and $f(r)=r(1+r^2)^{-1/4}$. 
Now using the Ricci curvature formulas 
(\ref{eq: ric-first})-(\ref{eq: ric-last}), we conclude that 
all the mixed terms vanish except $\Ric_g(Y_i,Y_j)$, and
$|\Ric_g(Y_i,Y_j)|\le ch^2$ when $i\neq j$.

Furthermore, by a straightforward computation
\begin{eqnarray*}
\Ric_g(Y_i,Y_i) & \geq & h^2\left(r^2(pK-L)+pR-S\right),  \\
\Ric_g(\pt_r,\pt_r) & = &  h^2\left(r^2(pK-L)+pR-S\right), \\ 
\Ric_g(U_i,U_i) & = & (p-2)f^{-2}(1-(f^\prime)^2)+h^2(Kr^2+R),  
\end{eqnarray*} 
where $K,L,R,S$
are some explicit linear functions of $n,c,m$ depending
on a particular element of $\beta_r$ (may be different in different equations 
above), and $K, R$ are positive. 
Note that $f^{-2}(1-(f^\prime)^2) \geq h^2 (\frac{3}{2} + r^2)$. 
Therefore, the matrix of the Ricci tensor in the basis $\beta_r$
is positive definite for all sufficiently large $p>k(n,c,m)$,
where $k(n,c,m)$ depends on $n, c, m$ and is independent of $r$. 
\end{proof}

\section{Metrics on fiber bundles}
\label{sec: metrics on fb}
This section provides an easy route to the proof 
of Corollary~\ref{cor: vect over ric nneg}.

\begin{thm} 
\label{thm: fb bounded error}
Let $\pi\co E\to B$ be a Riemannian submersion 
with totally geodesic fibers isometric to $F$, 
and a bounded error term. 
Assume that $E$ is complete, and
$B, F$ have nonnegative Ricci curvature.
Then $E\times \mathbb R^p$ admits a complete Riemannian metric of
positive Ricci curvature for all sufficiently large $p$.
\end{thm}
\begin{proof} Write $TE$ as the sum $\mathcal V\oplus\mathcal H$
of a vertical and a horizontal subbundles.
Fix a point $e\in E$.
By Section~\ref{sec: basis}, we can find an 
orthonormal basis $\{ X_k\}$ of vector fields 
on a neighborhood of $e\in E$ with the following
properties:\\
(1) each $X_k$ is either horizontal or vertical,\\
(2) if $X_k\in\mathcal H$, then $X_k$ is a horizontal lift of a vector field
$\check X_k$ on $B$,\\
(3) $\Ric_{B}(\check X_k, \check X_j)|_e=0$ for any $X_k, X_j\in \mathcal H$ 
with $k\neq j$,\\
(4) $\Ric_{F}(\hat X_k, \hat X_j)|_e=0$ for any  
$X_k, X_j\in\mathcal V$ with $k\neq j$,\\
(5) $\lp [X_k, X_j],X_k\rp|_e=0$ for all $k, j$.

For $t\in (0,1]$, define a metric on $E$ by $g_{t,E}$ where
$g_{t,E}(X+V,Y+U)=g_E(X,Y)+t^2 g_E(U,V)$ where
$g_E$ is the original metric on $E$, 
$X,Y\in \mathcal H$, and $U,V\in \mathcal V$.
Setting $Y_k=X_k$ if $X_k\in\mathcal H$, and 
$Y_k=X_k/t$ if $X_k\in\mathcal V$,
we get a $g_{t,E}$-orthonormal basis 
$\{Y_k\}$ on a neighborhood of $e$.   

If $C>0$ bounds the error term of $\pi$,
then the Ricci tensor $\Ric_{g_{t,E}}$ computed in the basis $\{Y_k\}$ 
satisfies $|\Ric_{g_{t,E}}(Y_k,Y_j)|\le Ct$
if $k\neq j$, and $\Ric_{g_{t,E}}(Y_k,Y_k)\ge -Ct^2 \geq -Ct$ for each $k$
(where we used $\Ric_B\ge 0$, $\Ric_F\ge 0$ for the latter
inequality). 
Now we are in position to apply Theorem~\ref{thm: positive} 
which completes the proof.
\end{proof}

\begin{cor}\label{cor: fb}
Let $B,F$ be compact Riemannian manifolds of nonnegative
Ricci curvature, and let $E\to B$ be a smooth fiber bundle
with fiber $F$ and structure group in the isometry group of $F$.
Then $E\times\mathbb R^p$ admits a complete Riemannian metric of
positive Ricci curvature for all sufficiently large $p$.
\end{cor}
\begin{proof}
By Example~\ref{ex: local metric} $E$ has a metric
making the projection $E\to B$ a Riemannian submersion with 
totally geodesic fibers isometric to $F$. Since $E$ is compact,
the submersion has a bounded error term, 
so we are done by Theorem~\ref{thm: fb bounded error}.
\end{proof}

\begin{proof}[Proof of Corollary~\ref{cor: vect over ric nneg}]
By Example~\ref{ex: vb error term} there is a complete Riemannian metric
on $E$ making the bundle
projection $E\to B$ a Riemannian submersion with totally geodesic
fibers and bounded error term so we are done by Theorem~\ref{thm: fb bounded 
error}.
\end{proof}

\begin{rmk}
It would be interesting to get a realistic lower bound on $p$.
The bound we get is an explicit function of $C$ and
$\dim(B)$, $\dim(F)$ which has little practical value
since we do not know how to estimate $C$. 
In the case (covered by Theorem~\ref{thm: vb over b}) when $B$ is a nilmanifold 
the constant $C$ can be estimated in terms of $\dim(B)$, and the structure
constants of the simply-connected nilpotent group covering $B$.
This indicates that $C$ cannot be chosen independently of 
the vector bundle.
\end{rmk}

\section{Metrics on iterated fiber bundles I}
\label{sec: ibi}
In this section we prove a technical lemma
which implies Theorem~\ref{thm: vb over b}
(in fact, in Theorem~\ref{thm: vb}
we deal with a larger class of base spaces). 
To prove Theorem~\ref{thm: vb over b}
we show that any manifold in $\mathcal B$ belongs to the
class $\mathcal M$ of manifolds of
almost nonnegative Ricci curvature 
with good local basis, so that Theorem~\ref{thm: positive} applies.

Let $G$ be a Lie group, $B$ be a smooth manifold, 
and $P\to B$ be a principal $G$-bundle over $B$
with a connection. 
Let $\pi\co E\to B$ be a smooth fiber bundle with fiber $F$ 
associated with $P$. The connection on $P$ defines
a decomposition of the tangent bundle to $E$ as the sum of
a vertical and a horizontal subbundles, $\mathcal V$ 
and $\mathcal H$.
Let $(E,g^E_t)$, $(B,g^B_t)$ be smooth families of complete 
Riemannian metrics parametrized by $t\in (0,1]$, and such that
for each $t$, $\pi$ is a Riemannian submersion with totally geodesic 
fibers and horizontal space $\mathcal H$. 
We denote the induced metric on $F$ by $g_t^F$. 

\begin{lem} \label{lem: main}
Let $\pi\co (E,g^E_t)\to (B,g^B_t)$ be a family
of Riemannian submersions as above.
Assume that $(B,g_t^B)\in\mathcal M_q(c_B, m_B)$, $(F,g_t^F)\in \mathcal 
M_r(c_F, m_F)$, 
and there are nonnegative constants
$L_a$, $L_b$, $L_f$, $b$, $f$ such that\\ 
(1) $|A_{X}Y|_{g_t^E}\le L_a$ for $t=1$ and all pairs of
unit orthogonal vector fields $X,Y\in\mathcal H$,\\
(2) $|K_{g_t^B}|\le L_b/t^{b}$ and $|K_{g_t^F}|\le L_f/t^{f}$.\\
Then the following holds for some positive constants $Q_4, Q_5$ and
$m=\max\{b, 2m_B, f\}$.\\
(A) If $r\ge 2m+3q$, then
$(E,\tilde g^E_t)\in\mathcal M_q(c_E,m_E)$ for some $\tilde g^E_t$ satisfying
$|K_{\tilde g^E_t}|\le Q_4/t^{2m+2q+f}$, where $m_E = \max \{m_B, m_F \}$.\\
(B) If $L_f=0=b$, then for each $p$,
$(E,\tilde g_{t}^{E,p})\in\mathcal M_p$ for some $\tilde g^{E,p}_{t}$ 
satisfying $|K_{\tilde g^{E,p}_{t}}|\le Q_5$. 
(C) If the $A$-tensor of $\pi$ is everywhere zero (i.e. $\pi$ is
a locally isometric flat bundle), then for each $p$,
$(E,\tilde g_{t}^{E,p})\in\mathcal M_p$ for some $\tilde g^{E,p}_{t}$ 
satisfying $|K_{\tilde g^{E,p}_{t}}|\le L_b t^{-pb/k}+L_f t^{-pf/k}$
where $k=\min\{r,q\}$. 
\end{lem}
\begin{proof}
For $s\in (0,1]$, define a metric on $E$ by 
$g_{t,s}^E(X+U,Y+V)=g_t^E(X,Y)+s^2g_t^E(U,V)$
where $X,Y\in \mathcal H$, and $U,V\in \mathcal V$.
Fix a point $e\in E$, and follow the procedure in
section~\ref{sec: basis} to 
construct an orthonormal basis $\{Y_i\}$ 
at a neighborhood of $e$ by combining  
the bases at $\pi(e)\in B$, and $e\in F$. In particular,
$\lp[Y_i,Y_j],Y_i\rp_{g_t}$ vanishes at $e$.

First, we estimate the $A$-tensor of $\pi$.
Since  $A_X Y=[X,Y]^{\mathcal V}$ for $X,Y\in\mathcal H$, and
the subbundles $\mathcal H$, $\mathcal V$ are independent of $t$, 
so is $A_X Y$.
For $Y_i, Y_j\in\mathcal H$ we get:
\[|A_{Y_i}Y_j|_{g_t^E}=t^{-m_i-m_j}|A_{X_i} X_j|_{g_t^E}\le 
t^{-m_i-m_j}|A_{X_i} X_j|_{g_1^E}\le L_a t^{-m_i-m_j}\le 
L_a t^{-2m_B},\]
where the first inequality uses that the $g_t$-length of any vector
at $e$ is bounded above by its $g_1$-length, since $t\le 1$.
Thus, $|A_X Y|_{g_t^E}\le (\dim B)^2 L_a t^{-2m_B}$,
for any unit vector fields $X,Y\in\mathcal H$.
Now let $U\in\mathcal V$. 
Since $A_{X} U$ is horizontal,
and $\lp A_{X} U, Y_j\rp_{g_t^E}=-\lp A_{X} Y_j, U\rp_{g_t^E}$, 
we deduce by computing in the basis $\{Y_j\}$ that 
$|A_{X} U|_{g_t^E}\le (\dim B)^2 L_a t^{-2m_B}$.

From the O'Neill formulas for the sectional 
curvature~\cite[9.29]{Bes} of $g_{t,s}^E$
we now see that
\[|K_{g_{t,s}^E}|\le 
|K_{g_t^B}|+\frac{4(\dim B)^2 L_a s^2}{t^{2m_B}}+
\frac{|K_{g_t^F}|}{s^2}\le 
\frac{L_b}{t^{b}}+\frac{4(\dim B)^2 L_a s^2}{t^{2m_B}}+
\frac{L_f}{s^2t^{f}}=L(t,s).\]
Each component of the Ricci tensor $\Ric_{g_{t,s}^E}(Y_i,Y_j)$
is the sum the components of the curvature tensor, which in turn
is the sum of sectional curvatures.
Therefore, we get the estimate:  
$|\Ric_{g_{t,s}^E}(Y_i,Y_j)|\le Q_1\cdot L(t,s)$ 
for any $i,j$ and some constant $Q_1> 0$.
Setting $s=1$, and using $|\Ric_{g_t^B}|\le (\dim B)L_b/t^b$ and
$|\Ric_{g_t^F}|\le (\dim F)L_f/t^f$, we conclude that
the error of the Riemannian submersion
$\pi\co (E,g_{t,1}^E)\to (B,g_{t,1}^B)$ is bounded above
by $Q_2\cdot L(t,1)$ for some constant $Q_2>0$.

To prove (A), let $m=\max\{b, 2m_B, f\}$ so that 
$Q_2\cdot L(t,1)\le Q_3/t^m$ for a constant $Q_3>0$. 
Setting $\tilde g_t^E=g^E_{t,s}$, where $s=t^{m+q}$,
we get from (\ref{eqn: 5.2})-(\ref{eqn: 5.4}) that
$|\Ric_{\tilde g_t^E}(Y_i,Y_j)|\le Q_3 t^q$ for any $i\neq j$,  
$\Ric_{\tilde g_t^E}(Y_i,Y_i)\ge -(c_B+Q_3)t^q$ for any $Y_i\in\mathcal H$, and
for any $Y_i\in\mathcal V$,
\[\Ric(Y_i,Y_i)_{\tilde g_t^E}\ge 
-c_Ft^rt^{-2m-2q}\ge -c_Ft^q,\] where
the last inequality follows from $r\ge 2m+3q$. 
Thus, $(E,\tilde g_t^E)\in\mathcal M_q$. Note that
$|K_{\tilde g^E_t}|\le L(t,t^{m+q})\le Q_4/t^{2m+2q+f}$ for some $Q_4>0$.

To prove (B) note that $b=0=L_f$ implies  that
$L(t,s)=L_b+4(\dim B)^2 L_a s^2/t^{2m_B}$. 
Let $\bar g_t^E=g^E_{t,s}$ with $s=t^{2m_B+1}$ so 
$|K_{\bar g^E_t}|\le Q_5$ for a constant $Q_5>0$.
Now the error of $\pi\co (E,g_{t,1}^E)\to (B,g_{t,1}^B)$ 
is bounded by $Q_6t^{-2m_B}$ for a constant $Q_6>0$, so 
$|\Ric_{\bar g_t^E}(Y_i,Y_j)|\le Q_6t$ if $i\neq j$,
$\Ric_{\bar g_t^E}(Y_i,Y_i)\ge -Q_6t-c_B t^q\ge -t^{k}(Q_6+c_B)$
where $k= \min \{1, q\} >0$.  
Thus,  $(E,\bar g_t^E)\in\mathcal M_{k}$ so 
reparametrizing $\tilde g_t^{E,p} = \bar g^E_{t^{p/k}}$, 
we have $(E,\tilde g_t^{E,p})\in\mathcal M_p$ for any $p>0$ and 
$|K_{\tilde g^{E,p}_t}|\le Q_5$. 

To prove (C) note that the metric $g_{t,1}^E$
satisfies $|K_{\tilde g^{E}_{t}}|\le L_b t^{-b}+L_f t^{-f}$
because $L_a=0$, and $(E,\bar g_t^E)\in\mathcal M_{k}$
where $k= \min \{r, q\} >0$. 
Reparametrizing $\tilde g_t^{E,p} = \bar g^E_{t^{p/k}}$,
we get the desired metric.   
\end{proof}

\begin{thm}\label{thm: vb}
Let $B$ be a compact manifold with a family of metrics $g_t^B$
satisfying $(B,g_t^{B})\in\mathcal M_q$ 
and $|K_{g_t^B}|\le L_b/t^b$ for some $q>0$, $L_b, b\ge 0$. 
If $\mathbb R^n\to E\to B$ 
is a vector bundle over $B$, then $E\in \mathcal M_q$.
\end{thm}
\begin{proof}
Fix $q>0$.
By Example~\ref{ex: vb error term} we make
the projection $E\to B$ into a Riemannian submersion.
Namely, consider the principal $O(n)$-bundle $P\to B$
such that $E\to B$ is associated with $P$, fix a connection
on $P$, and define a family of metrics on $P$ 
as in Example~\ref{ex: local metric}; in particular, the fiber is
$O(n)$ with an independent of $t$ biinvariant metric and
$P\to B$ is a Riemannian submersion with totally geodesic fibers
for each $t$.
Then $E$ gets the metric as the base of the Riemannian submersion
$(P\times\mathbb R^n)\to (P\times\mathbb R^n)/O(n)$.
Now the projection $E\to B$ is a Riemannian submersion with
totally geodesic fibers isometric to $(O(n)\times\mathbb R^n)/O(n)$.
The metric on the fibers is independent of $t$, and has sectional
curvatures within $[0,L_f]$ for some $L_f>0$.  
Note that the horizontal and vertical spaces on $P$
are independent of $t$, hence $E$ enjoys the same property.
As we mentioned in Example~\ref{ex: vb error term}, the $A$-tensor of
$E\to B$ satisfies Lemma~\ref{lem: main}(1).
It remains to apply Lemma~\ref{lem: main}(A).
\end{proof}

\begin{proof}[Proof of Theorem~\ref{thm: vb over b}]
By Theorems~\ref{thm: positive},~\ref{thm: vb}, it suffices to
find a family of metrics $g_t^B$ on $B$
satisfying $(B,g_t^{B})\in\mathcal M_q$ 
and $|K_{g_t^B}|\le L_b/t^b$. 
We argue inductively from the definition of the class $\mathcal B$.  
Given metrics on $F,B$, we always equip $E$ 
with the metric as in Example~\ref{ex: local metric}.
If $B$ is compact with $\Ric\ge 0$, then
$B\in \mathcal M_r$ for any $r$ by Example~\ref{ex: ric nneg}, and 
the sectional
curvature of $B$ is bounded by a constant. The induction step 
follows from Lemma~\ref{lem: main}(A), where the assumption
(1) holds by compactness. 
\end{proof}

\section{Metrics on iterated fiber bundles II}
\label{sec: ibii}
In this section we give more sophisticated examples of base
manifolds for which Theorem~\ref{thm: vb} applies. 
Roughly speaking, we allow base manifolds to be iterated bundles with almost 
nonnegatively curved fibers.
We first introduce some definitions which 
should make it easier to digest our results. 

Let $\mathcal S_{r,e}$ be the class
of smooth manifolds such that each $E\in\mathcal S_{r,e}$ admits
a smooth family of Riemannian metrics $g_t^E$ with 
$(E,g_t^E)\in\mathcal M_r$
and $|K_{g_t^{E}}|\le L_e/t^e$ for some
constant $L_e$. 
In these notations Theorem~\ref{thm: vb} says that
the total space of
any vector bundle over a compact manifold in $\mathcal S_{r,e}$ 
lies in $\mathcal M_r$.
By reparametrizing, one sees that
$\mathcal S_{r,0}=\mathcal S_{q,0}$ for any $r,q>0$,
so we denote $\mathcal S_{\infty,0}=
\cap_{r>0}\mathcal S_{r,0}=\cup_{r>0}\mathcal S_{r,0}$.
Thus if $E\in\mathcal S_{\infty,0}$, then for each $r$, 
$E$ has a family of metrics
$g_t^{E,r}$ with $(E,g_t^{E,r})\in \mathcal S_{r,0}$.

Given a collection $g_\alpha$ of Riemannian
metrics on a manifold $E$, we refer to the group 
$\cap_{\alpha}\Iso(E,g_\alpha)$ 
as the {\it symmetry group} of $(E,g_\alpha)$ and denote it
by $\mathrm{Sym}(E,g_\alpha)$.
In particular, if $(E,g^E_t)\in\mathcal S_{r,e}$, then
$\mathrm{Sym}(E,g^E_t)$ is $\cap_{t}\Iso(E,g_t^E)$, 
and if $(E,g^{E,r}_t)\in\mathcal S_{\infty,e}$, then
$\mathrm{Sym}(E,g^{E,r}_t)$ is $\cap_{t, r}\Iso(E,g_t^{E,r})$.

\begin{ex}\ \\
(1) If $\Ric(E,g)\ge 0$, then by Example~\ref{ex: ric nneg}, 
$\mathrm{Sym}(E,g)=\Iso(E,g)$ is the symmetry group of the 
family consisting of the metric $g$, and
$(E,g)\in\mathcal S_{\infty,0}$.\\
(2) Let $E$ be the total space of a
fiber bundle $F\to E\to B$ where $B$ is a compact manifold in
$\mathcal S_{\infty,0}$,
$F$ is a compact flat manifold, and the structure group of the bundle
lies in $\Iso(F)$. Then 
$E\in\mathcal S_{\infty,0}$ by Lemma~\ref{lem: main}(B)
where we equip $E$ with the metric as in 
Example~\ref{ex: local metric}.\\ 
(3) Let $E$ be the total space of a flat
fiber bundle $F\to E\to B$ where $(B,g_t^{B,r})\in\mathcal S_{\infty,0}$, 
$(F,g_t^{F,r})\in\mathcal S_{\infty,0}$,
and the structure group of the bundle lies in
$\mathrm{Sym}(F,g_t^{F,r})$. Equip $E$ with the metric as in 
Example~\ref{ex: flat}. 
If $B,F$ are compact, then
by Lemma~\ref{lem: main}(C),
$E\in\mathcal S_{\infty,0}$.\\
(4) Let $E$ be the total space of a
fiber bundle $F\to E\to B$ where $(B,g_t^{B})\in\mathcal S_{q,b}$, 
$(F,g_t^{F,r})\in\mathcal S_{\infty,0}$,
and the structure group of the bundle lies in
$\mathrm{Sym}(F,g_t^{F,r})$. If $B,F$ are compact, then
by Lemma~\ref{lem: main}(A),
$E\in\mathcal S_{p,e}$ for some $p,e$.\\
(5) Let $E$ be the total space of a flat
fiber bundle $F\to E\to B$ where $(B,g_t^{B})\in\mathcal S_{q,b}$, 
$(F,g_t^{F})\in\mathcal S_{r,f}$,
and the structure group of the bundle lies in
$\mathrm{Sym}(F,g_t^{F})$. Equip $E$ with the metric as in 
Example~\ref{ex: flat}. 
If $B,F$ are compact, then
by Lemma~\ref{lem: main}(C), $E\in\mathcal S_{p,e}$ for some $p,e$.
\end{ex}
 
\begin{rmk}
As one may suspect, the symmetry group is often discrete
so that any bundle with the structure group being the symmetry 
group of the fiber is almost the product. However,
here is an example with a large symmetry group.
Let $(B,g_t^B)\in\mathcal S_{r,e}$ and let $G$ be a compact
Lie group with a biinvariant metric. Then 
the total space $E$ of any principal $G$-bundle over $B$, 
equipped with the family of metrics $g^E_t$
as in Example~\ref{ex: local metric},
has the symmetry group containing $G$.
Furthermore, if $g^{E,r}_t$ is any reparametrization of $g_t^E$,
then of course $G\le\Iso(g^{E,r}_t)$.
\end{rmk}

\section{Fiber bundles with no nonnegative Ricci curvature}

\label{sec :fb and ric}
In this section we justify the claim made in the introduction
that if $E$ is a compact manifold with $\Ric(E)\ge 0$
which fibers over a torus $T$, then 
the pullback of the bundle $E\to T$ to a finite cover of $T$ 
has a section.
Our only tool is the fact that a finite cover of $E$ is homeomorphic to
the product of a simply-connected manifold and a torus~\cite{CG}.
In fact, we prove a more general result as follows.

Let $E$ be a manifold such that a finite cover 
$p\co \bar E\to E$ is homeomorphic to the product of a simply-connected 
manifold $C$, and a $k$-torus $T^k$. 
Fix an inclusion $i\co T^k\to \bar E$.
Assume that there is a torus $T$, and a homomorphism
$\phi\co\pi_1(E)\to\pi_1(T)$ with finite cokernel. 
(Such a homomorphism always exists if $E$ is a compact manifold
that fibers over a torus thanks to the homotopy exact
sequence of the fibration).

Since $T$ is aspherical, $\phi$ is induced by a 
continuous map $f\co E\to T$. 
The composition 
$\phi\circ p_*\circ i_*\co \pi_1(T^k)\to\pi_1(T)$
also has finite cokernel, so there is a finite cover
$\pi\co\tilde T\to T$ and a lift $q\co T^k\to \tilde T$
of $f\circ p\circ i$ such that $q$ is $\pi_1$-surjective.
Any surjection of free abelian groups has a section,
and since $T^k$ is aspherical, $q$ has a homotopy section.
Composing the homotopy section with $p\circ i$, we get a map 
$s\co\tilde T\to E$ such that $f\circ s$ is homotopic to $\pi$.

Replacing $E$ by a homotopy equivalent space $E^\prime$, we can
think of $f\co E\to T$ as a Serre fibration
$f^\prime\co E^\prime\to T$.
Then the pullback of the fibration via $\pi$ has a homotopy
section induced by $s$, and by the covering homotopy theorem
$f^\prime$ has a section.
Similarly, if $E\to T$ is a fiber bundle to begin with,
then its $\pi$-pullback has a section.

\begin{ex}\ \\
(1) Let $E\to T$ be a principal $G$-bundle over a torus $T$ 
where $G$ is a compact Lie group. If $\Ric(E)\ge 0$,
then $E\to T$ becomes trivial in a finite cover
because any principal bundle with a section is trivial. 
Using obstruction theory as in~\cite[Section 4]{BK2}, 
it is easy to find $G$-bundles
over $T$ which do not become trivial in a finite cover.\\
(2) Let $E\to T$ be a sphere bundle over a torus with
a nonzero rational Euler class. Then since finite covers
induce injective maps on rational cohomology, 
$E$ does not admit a metric with $\Ric\ge 0$ 
(see~\cite[Section 6]{BK2} for related results).
\\
(3) Let $E$ be an iterated fiber bundle, i.e.
$E=E_0\to E_1\to\dots\to E_n$ where 
$E_{i-1}\to E_i$ is a fiber bundle for each $i$. 
Assume that $E$ is a compact manifold, and $E_n=T$ is a torus. 
Since the composition of all the bundle projections $E\to T$
is a fiber bundle, it has a section if $\Ric(E)\ge 0$,
in which case the bundle $E_i\to T$ has a section for each $i$.
Thus, if $E_{n-1}\to E_n=T$ is a bundle
which does not have a section in a finite cover (e.g. as in (1) or (2)),
then $E$ admits no metric with $\Ric\ge 0$. 
\end{ex}

\small
\bibliographystyle{amsalpha}
\bibliography{base}

\

SCHOOL OF MATHEMATICS, GEORGIA INSTITUTE OF TECHNOLOGY,
ATLANTA, GA 30332-0160, USA

{\normalsize
{\it email:} \texttt{ib@math.gatech.edu}}

\

DEPARTMENT OF MATHEMATICS, UNIVERSITY OF CALIFORNIA SANTA BARBARA,

SANTA BARBARA, CA 93106, USA

{\normalsize
{\it email:} \texttt{wei@math.ucsb.edu}}
\end{document}